\documentclass[a4paper,12pt]{amsart}
\usepackage{amsmath,amsfonts,amssymb}
\usepackage{amsthm}
\usepackage{graphics,graphicx,epsf}
\usepackage[usenames]{color}
\usepackage{color}

\usepackage{tikz}
\usepackage{mathrsfs}

\usepackage{hyperref}

\usepackage[pagewise]{lineno}
\nolinenumbers
\newtheorem{theorem}{Theorem}[section]

\newtheorem{lemma}[theorem]{Lemma}
\newtheorem{cor}[theorem]{Corollary}
\newtheorem{rem}[theorem]{Remark}
\newtheorem{ex}[theorem]{Example}

\renewcommand{\proof}{{\noindent \bf Proof. \ }}

\newcommand{\eproof}{\hfill \mbox{${\square}$}}
\newcommand{\R}{\mathbb{R}}

\newcommand{\dis}{\displaystyle}
\newcommand{\dist}{\textrm{dist}}

\begin{document}

\title[\hfil  Neumann condition for a nonlocal  problem]
{A note on the formulation of the Neumann boundary condition for a nonlocal diffusion problem with continuous kernel}

\author[A. L. Pereira]
{Ant\^onio L. Pereira}

\address{Ant\^onio L. Pereira \hfill\break
Instituto de Matem\'atica e Estat\'istica,
Universidade de S\~ao Paulo,
S\~ao Paulo, Brazil}
\email{alpereir@ime.usp.br}

\thanks {Partially supported by FAPESP-SP Brazil grant 2020/14075-6.}

\subjclass{45J05, 45H05, 37L05{, 45K05}}
\keywords{Nonlocal diffusion, convolution kernel, population model}

\date{January 31, 2026}

\begin{abstract}
The nonlocal diffusion equation
\begin{equation*}
      \begin{cases} u_t(t,x) & = \int_{\Omega}   K(x,y) \,
       u(t,y) \,  d\,y
      -  \int_{\Omega}   K(y,x) \, u(t,x) \, d \, y \\
      u(0,x) & = u_0(x)
      \end{cases}
     \end{equation*}
    has been proposed as a model for some evolution processes with diffusion, including population models. However, in general, this formulation does not satisfy
$  \int_{\Omega}   K(y,x)  \, d \, y = 1$, as
 expected  from  its interpretation  as a probability density.
 In this note, we propose a modification of the kernel, based on the idea of `reflection' at the boundary, a technique familiar in one dimensional problems. We show that a similar construction is possible in higher dimensions, with the new kernel satisfying the aforementioned integral equality and being  also symmetric in certain special cases.
\end{abstract}

\maketitle

\section{Introduction}

   The  evolution equation
   \begin{equation}
   \label{evolution1}
 \begin{cases} u_t(t,x) & = \int_{\R^n}   K(x,y) \, u(t,y) \,  d\,y
      - \int_{\R^n}  K(y,x) \,  u(t,x) \,  d\,y \\
      u(0,x) & = u_0(x),
      \end{cases}
      \end{equation}
   where  $K: \R^n \times \R^n  \to {[0,\infty)}$  is a regular integrable kernel,  has been proposed as a nonlocal analogue of the diffusion equation
    \begin{equation*}
      \begin{cases} u_t(t,x) & = \Delta u(t,x), \quad x\in \R^n \\
      u(0,x) & = u_0(x)
      \end{cases}
     \end{equation*}
  (see \cite{fife}).

    In population models, $K(x,y)$  can be interpreted as the probability density  for  a species moving from site $y$ to site $x$, which makes it reasonable to assume the  hypothesis
     \begin{equation} \label{outflux} \int_{\R^n} K(x,y) \, d\,x = 1
     , \quad \textrm{for any }  \ y \in \R^n,
     \end{equation}
     since this gives the probability of going from $y$ to an arbitrary location.
        In general, we do not  expect
     \begin{equation} \label{influx}
    \int_{\R^n} K(x,y) \, d\,y = 1,
      \end{equation}
        since this measures the influx at the site $x$, from all other points, and this  may vary (i.e., some sites may be `more attractive' than others).
   If { hypothesis  \eqref{outflux}} is assumed,  equation \eqref{evolution1} becomes
     \begin{equation*}
      \begin{cases} u_t(t,x) & = \int_{\R^n}   K(x,y) \,
       u(t,y) \,  d\,y
      -  u(t,x)  \\
      u(0,x) & = u_0(x),
      \end{cases}
     \end{equation*}
   which has been considered by many authors (see \cite{rossi} and references therein).

  Now, if we consider the problem in a bounded domain $\Omega $, some `boundary condition' must be added for the problem to be well posed.
   If we assume a hostile environment, the species dies if it moves to a site outside $\Omega$. We then impose the solutions to be
    identically zero outside the domain $\Omega$, obtaining a nonlocal version of the Dirichlet problem. The evolution equation is the same,  except that  the integration is restricted to the domain $\Omega$, that is
  \begin{equation*}
      \begin{cases} u_t(t,x) & = \int_{\Omega}   K(x,y) \,
       u(t,y) \,  d\,y
      -  u(t,x), \textrm{ if }  x \in \Omega,  \\
      u(t,x) & =0,  \textrm{ if }  x \in \R^n \setminus \Omega, \\
      u(0,x) & = u_0(x) .
      \end{cases}
     \end{equation*}

   If we impose that  interactions occur only inside the domain, we should obtain a version of the Neumann problem. It is, however, not clear how exactly to formulate the problem in this case.  A possible choice (see, for example \cite{bernal} and  \cite{rossi}) is to restrict both integrals in \eqref{evolution1} to the domain $\Omega$, leading to the problem
   \begin{equation*}
      \begin{cases} u_t(t,x) & = \int_{\Omega}   K(x,y) \,
       u(t,y) \,  d\,y
      -  \int_{\Omega}   K(y,x) \, u(t,x) \, d \, y \\
      u(0,x) & = u_0(x) .
      \end{cases}
     \end{equation*}

      This formulation has some nice features. In particular, if
     $ \displaystyle  \int_{\Omega}   K(x,y)  \,  d\,y
      =  \int_{\Omega}   K(y,x) \, d \, y  $, the constants are equilibrium solutions, which also happens in the local diffusion Neumann  problem. This hypothesis obviously holds  if $K(x,y)$ is symmetric and, in particular, if  $K(x,y) = J(x-y)$ is of convolution type and $J$ is even.

What seems somewhat unsatisfactory  with this formulation is that the integral
$ \dis \int_{\Omega}   K(y,x)  \, d \, y $ cannot {be} expected to be  identically equal to $1$, as
 required by  its interpretation  as the probability of transition from $x$ to any point in the domain, since the transition to points outside the domain is now forbidden.
 {In probabilistic terms, jumps that would leave $\Omega$ are simply suppressed (the individual remains at $x$). An alternative, which we explore here, is to redirect such jumps back into $\Omega$, in the spirit of a reflection at the boundary.}

   We may assume that the kernel $ \dis    K(x,y)   $, instead of  condition \eqref{outflux}, now satisfies
   \begin{equation} \label{outflux_bounded} \int_{\Omega} K(y,x) \, d\,y = 1, \quad \textrm{for any }  \ x \in \Omega.
     \end{equation}
{Under \eqref{outflux_bounded} the Neumann problem above reads $u_t(t,x) = \int_\Omega K(x,y)\,u(t,y)\,dy - u(t,x)$ and, by Fubini's theorem,
\[
\frac{d}{dt}\int_\Omega u(t,x)\,dx = \int_\Omega u(t,y)\Big(\int_\Omega K(x,y)\,dx\Big)dy - \int_\Omega u(t,x)\,dx = 0,
\]
so the total mass is conserved, as in the local Neumann problem; no symmetry of $K$ is needed for this.}

   In fact, starting with a kernel $ \dis    K(x,y)   $ satisfying    \eqref{outflux} we can
    modify it in various ways  in order to obtain a new {kernel}
    $ \dis    {\tilde{K}(x,y)}   $ satisfying \eqref{outflux_bounded}.  For instance, we may  define
      \begin{equation}
       \label{modify_kernel1}
      \tilde{K}(y,x)  : =  \frac{1}{h(x)} K(y,x)  \textrm{ with }  h(x) :=
       \int_{\Omega} K(y,x) \, d\, y,
     \end{equation}
     which is well defined if $ h(x) \neq 0$.
     However, besides the artificiality and lack of motivation for this definition,  the kernel  thus defined generally  does not  satisfy
          $ \displaystyle  \int_{\Omega}   \tilde{K}(x,y)  \,  d\,y
      = 1 $, so the constants are not equilibria.  Also,  $\dis \tilde{K}(x,y) $ is not   symmetric, even if
  $\dis {K}(x,y) $  is symmetric. Ideally, one would like the kernel
  to satisfy, besides the hypothesis \eqref{outflux_bounded}, also the following  identity:
  \begin{align}
    \int_{\Omega} K(x,y) \, d\,y & =  \int_{\Omega} K(y,x) \, d\,y , \quad \textrm{for any }  \ x \in \Omega. \label{hyp_2}
  \end{align}
  This hypothesis is obviously fulfilled if $K$ is symmetric.

   In the one dimensional case, we  can construct a  kernel satisfying \eqref{outflux_bounded} and \eqref{hyp_2} by `reflecting at the boundaries'
    (see \cite{hutson} and  \cite{belletini}).
  Starting with a kernel  of the form $K(x,y) = J(x-y)$, where $J $ is even, has compact support in the  interval $[-T,T]$ { with $T>0$}, and
      $ \dis \int_{\R}J (z) \, d\,z =1$,
        we define  the Neumann kernel in the interval{ $ I=]0,T[$ } by
\[ K^{N}(x,y)  :=  J(x-y) +   J(x-R(y)) +  J(x- L (y)) \]
   where $ R(y)= T +( T-y) = 2 T-y $  is the reflection of $y$  with respect to the right end of $I$   and
  $ L(y)= -y$    is the reflection of $y$  with respect to the left end of $I$.
  {Since $J$ is even, $J(x-R(y)) = J(R(x)-y)$ and $J(x-L(y)) = J(L(x)-y)$, so that
  $K^N(x,y) = J(x-y) + J(R(x)-y) + J(L(x)-y)$, which is the form adopted in Section~2.}

   This construction is more natural than the one proposed in
   \eqref{modify_kernel1} above and, more importantly, the   kernel $K^N$
     is still symmetric
     and satisfies
 \begin{equation*}  \int_{0}^T K^N (x,y) \, d\,y = {\int_{x-2T}^{x+T} J (z) \, d\,z =} \int_{-T}^T J (z) \, d\,z =  1
    \end{equation*}
    for any $x \in [0,T]$, so \eqref{outflux_bounded} and \eqref{hyp_2} are met.

  The purpose of this note is to  show that   a similar construction is possible in higher dimensions, with the   resulting kernel satisfying condition \eqref{outflux_bounded}  and also \eqref{hyp_2} in the case of some special domains, when  it is also symmetric.

\section{Reflecting kernel for $\mathcal{C}^2$ domains}

      Suppose $\Omega \subset \R^n$ is a bounded domain, and the {kernel $K(x,y) $ } {satisfies \eqref{outflux} and} vanishes outside a neighborhood of the diagonal. More precisely  $K(x,y) = 0$ if $\| x-y \| {\geq} \delta $, for some $\delta>0$.
 {If $U_1, \dots, U_m$ are measurable subsets of $\Omega$ and $\varphi_k : U_k \to \R^n$, $k=1,\dots,m$, are measurable maps,}
we may define a new kernel $\tilde{K}(x,y)$ by
     \begin{equation} \label{modf_kernel1}
     \tilde{K}(x,y) = K(x,y) + {\sum_{k=1}^{m} \chi_{U_k}(x)\, K(\varphi_{k}(x),y)},
\end{equation}
{where $\chi_{U_k}$ is the characteristic function of $U_k$ (the $k$-th term is taken to be zero for $x \notin U_k$).}

We now investigate some properties of $ \tilde{K}$, under appropriate conditions on  the maps  $\varphi_k$.

  \begin{lemma}\label{const_integral}
  {Let $U_1, \dots, U_m \subset \Omega$ and $\varphi_k: U_k \to \R^n$ be as above, and suppose that:}
  \begin{enumerate}
  \item {$\varphi_k$ is injective and $\varphi_k(U_k) \subset \Omega^c$, for $k=1, \dots, m$, and the sets $\varphi_k(U_k)$ are pairwise disjoint;}
   \item {the set $\{ z \in \Omega^c : \dist(z, \Omega) < \delta\}$ is contained in $\bigcup_{k=1}^m \varphi_k(U_k)$, up to a set of Lebesgue measure zero;}
   \item {each $\varphi_k$ preserves Lebesgue measure, that is,
   $\int_{U_k} f(\varphi_k(x))\, dx = \int_{\varphi_k(U_k)} f(z)\, dz$ for every measurable $f \geq 0$.
   (This holds, for instance, if $\varphi_k$ is a $\mathcal{C}^1$ diffeomorphism of the open set $U_k$ onto its image with $|\det D\varphi_k| \equiv 1$.)}
  \end{enumerate}
  Then
  \[ \int_{\Omega}   \tilde{K}(x,y)\, d \,x    =  \int_{\R^n}   K(x,y)\, d \,x = 1, \textrm {\ for any } y\in \Omega.\]
 \end{lemma}

  \proof
  {Fix $y \in \Omega$. Then
   \begin{align*}
         \int_{\Omega}   \tilde{K}(x,y)\, d \,x   & =
          \int_{\Omega}   {K}(x,y)\, d \,x  + \sum_{k=1}^m
           \int_{U_k}   {K}(\varphi_k(x),y)\, d \,x \\
               & =       \int_{\Omega}    K(x,y) \, d \, x  +
 \sum_{k=1}^m   \int_{\varphi_{k}(U_k)} K(z,y) \, dz   \\
  & =   \int_{\Omega}    K(x,y) \, d \, x  +   \int_{ \bigcup_k \varphi_k (U_k) } K(z,y) \, dz  \\
    & =  \int_{\Omega}    K(x,y) \, d \, x  + \int_{ \Omega^c } K(z,y) \, dz  \\
    & =  \int_{\R^n}    K(x,y) \, d \, x = 1.
  \end{align*}
  We have used hypothesis 3 in the second line and the disjointness in hypothesis 1 in the third line.
  In the fourth line we used that $\bigcup_k \varphi_k(U_k) \subset \Omega^c$ (hypothesis 1), that
  $K(z,y) = 0$ if $\|z-y\| \geq \delta$, so that only points of $\Omega^c$ with $\dist(z,\Omega) < \delta$ contribute to either integral, and hypothesis 2. The last line is \eqref{outflux}.}
 \eproof

     \begin{rem}
     Under the conditions of Lemma \ref{const_integral},  the points ${\varphi_k}(x)$ can be interpreted  as a  {forbidden}  location in $\Omega^c$ for  {an individual at} $y$, which then  goes to
     $x$ instead. {In particular, $\tilde K$ satisfies \eqref{outflux_bounded} (with the roles of the variables as in \eqref{outflux}), so the Neumann problem with kernel $\tilde K$ reads $u_t = \int_\Omega \tilde K(x,y)\, u(t,y)\, dy - u(t,x)$ and conserves mass.}
     \end{rem}

 The natural question now is whether one can find reasonable assumptions  on the domain $\Omega$ and the kernel $K(x,y)$ so that the conditions of { Lemma \ref{const_integral} } are fulfilled.

  We consider first a class of  domains in $\R^2$.

   \begin{ex} \label{ex_starshaped}
   {Let $\rho: \R \to (0,\infty)$ be a $2\pi$-periodic $\mathcal{C}^1$ function and let
   $R \subset \R^2$ be the (star-shaped) domain given in polar coordinates by}
   $$ R :=  {\{ P(r, \theta) \ : \   0 \leq  r  <  \rho(\theta),\ 0  \leq  \theta < 2 \pi\}}, $$
    where $ P(r, \theta) = (r \cos \theta, r \sin \theta)$ is the polar transformation of coordinates.
    {We look for a map of the form $\Psi(r,\theta) = (g_\theta(r), \theta)$ and
    $\Phi := P \circ \Psi \circ P^{-1}$ on $R \setminus\{0\}$.
    Since $JP(r,\theta) = r$ and $D\Psi$ is triangular with $\det D\Psi(r,\theta) = g_\theta'(r)$, we obtain
    \[ \det D\Phi(z) = \frac{g_\theta(r)\, g_\theta'(r)}{r}, \qquad (r,\theta) = P^{-1}(z).\]
    Imposing $\det D\Phi \equiv -1$ and $g_\theta(\rho(\theta)) = \rho(\theta)$ (points of $\partial R$ are fixed), we arrive at
    $\big(g_\theta^2\big)' = -2r$, $g_\theta(\rho(\theta)) = \rho(\theta)$, whose solution is
    $g_\theta(r) = \sqrt{2\rho(\theta)^2 - r^2}$. Thus
    \[ \Phi(z) = \sqrt{2 \rho(\theta)^2 - \|z\|^2}\; \frac{z}{\|z\|}, \qquad z = P(r,\theta) \neq 0 .\]}

  The map $\Phi$ is a kind of reflection of each point {$P(r,\theta) \in R$} with respect to the point {$P(\rho(\theta),\theta)$}, composed with a squeezing in the radial direction in order to make it area preserving.
   {It maps $R \setminus \{0\}$ bijectively onto the region outside $R$
     $$ A^{\prime} = \{P(r, \theta) :  \rho(\theta) <  r <  \sqrt{2}\,   \rho(\theta) , \ 0 \leq \theta < 2 \pi \}, $$
 and, since $g_\theta(g_\theta(r)) = r$, it satisfies $\Phi \circ \Phi = \mathrm{id}$ on $(R\setminus\{0\}) \cup A'$.}
 {Let $\Gamma' = \{P(\sqrt 2\rho(\theta),\theta) : 0 \leq \theta < 2\pi\}$ be the outer boundary of $A'$ and
 $d := \dist(\partial R, \Gamma') > 0$. If $z \notin \overline{R} \cup A'$, any segment joining $z$ to a point of $\overline R$ crosses $\Gamma'$, so $\dist(z, R) \geq d$. Hence,
  if $0 < \delta \leq d$, the hypotheses of Lemma \ref{const_integral} are met with  $m=1$, $U_1 = R \setminus\{0\}$ and  $\varphi_1 = \Phi$ (note that $\partial R$ has measure zero).
  We remark that $d \leq (\sqrt2 - 1)\rho_{\min}$, but $d$ may be much smaller than $(\sqrt2-1)\rho_{\min}$ if $\partial R$ is very oblique with respect to the rays through the origin.} \eproof
    \end{ex}

 To extend this construction to more general domains, we  need the following result,  concerning the existence of
  `normal coordinates'.

\begin{theorem} \label{normal_neighborhood}  Let $ \Omega \subset \R^n$ be a {bounded} domain with a $\mathcal{C}^m$ boundary, $ m \geq 2$. There exists $r>0$ so that, if
  \begin{itemize}
  \item $B_r(\partial \Omega)= \{ x \in \R^n : \dist(x, \partial \Omega) < r  \}$,
 \item  $\pi(x) = \textrm{ the point of  }  \partial \Omega
 \textrm{ nearest to  } x  $,
 \item  $  t(x) =  \pm  \dist(x, \partial \Omega)$
  ($+$ outside, $-$ inside),
  \end{itemize}
then:
\begin{itemize}
\item  $ t(\cdot) : B_r(\partial \Omega) \to (-r,r) $ and
 $\pi(\cdot) :  B_r(\partial \Omega) \to \partial \Omega $ are
 well-defined,
 $\pi$ is  a $\mathcal{C}^{m-1} $ retraction onto
  $\partial \Omega$
and $t$ is of class $\mathcal{C}^{m} $;
 \item
$ x \mapsto  (\pi (x), t(x) ) :  B_r(\partial \Omega) \to \partial \Omega \times (-r,r)$  is a  $\mathcal{C}^{m-1} $
 diffeomorphism with inverse
$ \Psi(\xi, t) = \xi + t N(\xi)$,
 where $N(\xi) $ is the unique outward unit normal to $\partial \Omega$
  at $\xi$.
\end{itemize}
 \noindent Furthermore,
 $t(\cdot)$ is the unique solution of
 $| \nabla t(x)| = 1$ in   $B_r(\partial \Omega)$ with $t=0$ on
  $\partial \Omega$ and  $\frac{\partial t }{\partial  N } > 0 $ on  $\partial \Omega$.
  Extending the normal field $N$ to $B_r(\partial \Omega)$ by
$ N(\xi + t N(\xi) )= N(\xi)$, $-r <t < r $, we have
 $ N (x) = \nabla t(x)$  on $B_r (\partial \Omega).$ Also
 $ {\mathcal{S}} (x) = D N (x) = D^2 t(x)$, restricted to the tangent space at $ x  \in \partial \Omega$,  is the curvature {(shape operator)} of $\partial \Omega$.
\end{theorem}

\proof See \cite{henry}. \eproof

\begin{theorem} \label{area_preserving} Let $ \Omega \subset \R^n$ be a {bounded} domain with a $\mathcal{C}^m$ boundary, $ m \geq 2$,
{and let $r>0$ be as in Theorem \ref{normal_neighborhood}. Then there exist $0 < \varepsilon < r$, $\varepsilon' > 0$ and a map
$\Phi : U_\varepsilon \to \R^n$, where $U_\varepsilon := \{x \in \Omega : \dist(x,\partial\Omega) < \varepsilon\}$,
satisfying the following properties:}
 \begin{enumerate}
  \item {$\Phi$ is a homeomorphism of $U_\varepsilon$ onto an open set $V \subset \Omega^c$}{, where $\Omega^c$ is the complement set of $\Omega$ in $\R^n$}, {and $V \supset \{ z \in \R^n \setminus \overline\Omega : \dist(z,\partial\Omega) < \varepsilon'\}$;}
  \item {$\Phi$ preserves Lebesgue measure: $\int_{U_\varepsilon} f(\Phi(x))\,dx = \int_V f(z)\,dz$ for every measurable $f \geq 0$;}
  \item {$\Phi$ is the restriction of an involution: the formula below defines $\Phi$ also on $V$, and $\Phi\circ\Phi = \mathrm{id}$ on $U_\varepsilon \cup V$;}
  \item {if $m \geq 3$, $\Phi$ is a $\mathcal{C}^{m-2}$ diffeomorphism of $U_\varepsilon$ onto $V$ and $\det D\Phi (x) = -1$ for all $x \in U_\varepsilon$.}
  \end{enumerate}
\end{theorem}

\noindent \proof
{Let $\Psi(y,s) = y + sN(y)$, $(y,s) \in \partial\Omega\times(-r,r)$, be the inverse of the normal coordinates given by Theorem \ref{normal_neighborhood}.
Let $\lambda_1(y), \dots, \lambda_{n-1}(y)$ be the eigenvalues of the curvature $\mathcal{S}(y) = DN(y)$ restricted to $T_y\partial\Omega$, and let
$v_1, \dots, v_{n-1}$ be a corresponding orthonormal basis of eigenvectors of $T_y\partial\Omega$.
Since $D\Psi_{(y,s)}(v_i,0) = v_i + s\,\mathcal{S}(y) v_i = (1+s\lambda_i)\, v_i$ and $D\Psi_{(y,s)}(0,1) = N(y)$, we have
\[
 \det D\Psi_{(y,s)} = \det\big[(1+s\lambda_1)v_1, \dots, (1+s\lambda_{n-1})v_{n-1}, N(y)\big] = p_y(s),
\]
where
\[ p_y(s) := \prod_{i=1}^{n-1}(1+s\lambda_i(y)) = \det\big(I + s\,\mathcal{S}(y)\big|_{T_y\partial\Omega}\big),\]
the determinant being computed with respect to the product of the surface measure $dS$ on $\partial\Omega$ and $ds$. Since $\Psi$ is a diffeomorphism, $p_y(s) \neq 0$ for $|s| < r$, and $p_y(0) = 1$, so $p_y(s) > 0$. Thus, for every measurable $f \geq 0$ supported in $B_r(\partial\Omega)$,
\begin{equation}\label{normal_coord_integral}
 \int_{\R^n} f(x)\, dx = \int_{\partial\Omega}\int_{-r}^{r} f\big(y + sN(y)\big)\, p_y(s)\, ds\, dS(y).
\end{equation}
The coefficients of the polynomial $p_y$ are the elementary symmetric functions of the $\lambda_i(y)$, which are continuous (indeed $\mathcal{C}^{m-2}$) in $y$.}

{For each $y \in \partial\Omega$ consider the o.d.e.}
\begin{equation}
 \label{diff_equation} \begin{cases}
 \dfrac{d\varphi}{ds} = -\dfrac{(1+ s \lambda_1(y))(1+ s \lambda_2(y))
 \cdots(1+ s \lambda_{n-1}(y))}{(1+ \varphi \lambda_1(y))(1+ \varphi \lambda_2(y))
 \cdots(1+ \varphi \lambda_{n-1}(y))} {= -\dfrac{p_y(s)}{p_y(\varphi)}},\\[2mm]
 \varphi(y,0) = 0.
\end{cases} \end{equation}
{This equation is separable: writing $P_y(s) := \int_0^s p_y(\sigma)\, d\sigma$, which is a strictly increasing function on $(-r,r)$, its solution is
\[ \varphi(y,s) = P_y^{-1}\big(-P_y(s)\big). \]
By continuity and compactness of $\partial\Omega$, there are $c, C > 0$ with $P_y(r/2) \geq c$ and $0 < p_y(s) \leq C$ for all $y \in \partial\Omega$, $|s| \leq r/2$; choosing $0 < \varepsilon < \min\{r/2, c/C\}$, we have $|P_y(s)| < c$ for $|s| < \varepsilon$, so that $\varphi(y,s)$ is well defined, continuous in $(y,s)$, $\mathcal{C}^1$ in $s$, and $\varphi(y,s) \in (-r/2,r/2)$ for $|s| < \varepsilon$.
Moreover, $\partial_s\varphi < 0$, so $\varphi(y,\cdot)$ maps $(-\varepsilon, 0)$ decreasingly onto $(0, \varphi(y,-\varepsilon))$, and $P_y(\varphi(y,\varphi(y,s))) = -P_y(\varphi(y,s)) = P_y(s)$, that is, $\varphi(y,\cdot)$ is an involution.
Note that $P_y(s)$ is the volume, per unit area of $\partial\Omega$ at $y$, of the layer between $\partial \Omega$ and $\partial\Omega + sN$; thus $\varphi(y,s)$ is the height of the outer layer with the same volume as the inner layer of depth $|s|$.}

{Define $\xi(y,s) := (y, \varphi(y,s))$ and
  \begin{equation*}
   \Phi(x) := \Psi\circ\xi\circ\Psi^{-1}(x) = \pi(x) + \varphi( \pi(x), t(x))\, N(\pi(x)), \qquad x \in U_\varepsilon,
  \end{equation*}
and let $V := \Psi\big(\{(y,\sigma) : y \in\partial\Omega,\ 0 < \sigma < \varphi(y,-\varepsilon)\}\big)$.
Since $\xi$ is a homeomorphism of $\partial\Omega\times(-\varepsilon,0)$ onto $\{(y,\sigma): 0<\sigma<\varphi(y,-\varepsilon)\}$ (with inverse given by $\xi$ itself), $\Phi$ is a homeomorphism of $U_\varepsilon$ onto $V$, and $\Phi\circ\Phi = \mathrm{id}$ on $U_\varepsilon\cup V$. Since $t > 0$ on $V$, $V \subset \Omega^c$, and $V \supset \{z : 0 < t(z) < \varepsilon'\}$ with $\varepsilon' := \min_{y\in\partial\Omega} \varphi(y,-\varepsilon) > 0$. This proves 1 and 3.}

{To prove 2, let $f \geq 0$ be measurable. By \eqref{normal_coord_integral} and the change of variables $\sigma = \varphi(y,s)$ for fixed $y$, for which $p_y(s)\,ds = -p_y(\sigma)\,d\sigma$ by \eqref{diff_equation},
\begin{align*}
\int_{U_\varepsilon} f(\Phi(x))\, dx
 & = \int_{\partial\Omega}\int_{-\varepsilon}^{0} f\big(y + \varphi(y,s)N(y)\big)\, p_y(s)\, ds\, dS(y) \\
 & = \int_{\partial\Omega}\int_{0}^{\varphi(y,-\varepsilon)} f\big(y + \sigma N(y)\big)\, p_y(\sigma)\, d\sigma\, dS(y)
 = \int_V f(z)\, dz .
\end{align*}
Finally, if $m \geq 3$, then $\varphi$ is $\mathcal{C}^{m-2}$ in $(y,s)$ and so is $\Phi$. The matrix of $D\xi(y,s)$ is block triangular,
$D\xi(y,s)(v,0) = (v, D_y\varphi(y,s)\,v)$, $D\xi(y,s)(0,1) = (0, \partial_s\varphi(y,s))$, so $\det D\xi = \partial_s\varphi$, and
\[
 \det D\Phi(x) = \det D\Psi_{\xi(y,s)}\cdot \partial_s\varphi(y,s)\cdot \big(\det D\Psi_{(y,s)}\big)^{-1}
 = p_y(\varphi)\cdot\Big(-\frac{p_y(s)}{p_y(\varphi)}\Big)\cdot \frac{1}{p_y(s)} = -1,
\]
where $(y,s) = \Psi^{-1}(x)$ and $\varphi = \varphi(y,s)$. This proves 4.}
  \eproof

  \begin{cor}\label{suitable_kernel}
   Let $\Omega \subset \R^n$ be a {bounded} domain with a $\mathcal{C}^2$ boundary {and let $K$ satisfy \eqref{outflux} and}
 $K(x,y) = 0 $ if $ \|x-y\| {\geq} \delta $. {Let $\varepsilon$, $\varepsilon'$, $U_\varepsilon$ and $\Phi$ be given by Theorem \ref{area_preserving}.
 If $0 < \delta \leq \varepsilon'$, the kernel
 $ \tilde{K} (x,y) =  K(x,y) + \chi_{U_\varepsilon}(x)\, K (\Phi(x),y)$}
 satisfies
 \[  \int_{\Omega} \, \tilde{K}(x,y) \, d \,x =1, \textrm{  for any }  y \in \Omega.\]
  \end{cor}
  \proof {We apply Lemma \ref{const_integral} with $m = 1$, $U_1 = U_\varepsilon$ and $\varphi_1 = \Phi$. Hypotheses 1 and 3 follow from items 1 and 2 of Theorem \ref{area_preserving}. As for hypothesis 2, if $z \in \Omega^c \setminus \partial\Omega$ and $\dist(z,\Omega) < \delta \leq \varepsilon'$, then $0 < t(z) < \varepsilon'$, so $z \in V = \Phi(U_\varepsilon)$; since $\partial\Omega$ has measure zero, hypothesis 2 holds.}  \eproof

 The map $\Phi$ can be given more {explicitly} in concrete examples. We consider a particularly simple example.

   \begin{ex} \label{ex_ball}
   Let $\Omega = B_{\rho} \subset \R^n$ be the ball of radius $\rho$ in $\R^n$.
   Then $ \pi(x) = \rho \frac{x}{\|x\|} $, $t(x) = \|x\| - \rho $
    and $N(y) = \frac{y}{\rho}$ is the unit exterior normal at the point
 $y \in \partial \Omega$. {All principal curvatures are equal to $1/\rho$, and}
 the equation \eqref{diff_equation} now becomes
 \[
  \begin{cases}
 \dfrac{d\varphi}{ds} = -\dfrac{(1+ \frac{s}{{\rho}})^{n-1}}{(1+ \frac{\varphi}{{\rho}})^{n-1}},\\[2mm]
 \varphi(0) = 0,
\end{cases} \]
whose solution is given by
$$\varphi(s) = - \rho + \sqrt[n]{2\rho^n - (\rho+s)^n}, \quad
 {-\rho < s < (\sqrt[n]{2}-1)\, \rho}.$$
The function $\Phi$ is then defined for any $x \in \Omega\setminus \{0\} $, and given by
 \[\Phi(x) = \rho \frac{x}{\|x\|} + \big( -\rho + \sqrt[n]{2\rho^n -  \|x\|^n}  \big)
  \frac{x}{\|x\|} { = \sqrt[n]{2\rho^n - \|x\|^n}\, \frac{x}{\|x\|}}. \]
  The image of $\Phi$ is the annulus $ \{ x\in \R^n: \rho < \|x\| < \sqrt[n]{2} \rho \}$, so, in this case, the conclusions of Corollary \ref{suitable_kernel} hold {(applying Lemma \ref{const_integral} directly with $U_1 = \Omega\setminus\{0\}$)}, if
  $0< \delta \leq  (\sqrt[n]{2} -1) \rho $.  \eproof
   \end{ex}

\begin{rem}
{Examples \ref{ex_starshaped} and \ref{ex_ball} are particular cases of the following construction, which does not require $\mathcal{C}^2$ regularity. If $\Omega = \{ r\omega : \omega \in S^{n-1},\ 0 \leq r < \rho(\omega)\}$ is star-shaped with respect to the origin, with $\rho : S^{n-1}\to(0,\infty)$ continuous, then
\[ \Phi(r\omega) := \big(2\rho(\omega)^n - r^n\big)^{1/n}\, \omega \]
maps $\Omega\setminus\{0\}$ onto $\{r\omega : \rho(\omega) < r < 2^{1/n}\rho(\omega)\}$, is an involution and preserves Lebesgue measure, since the volume element is $r^{n-1}dr\, d\omega$ and $s = (2\rho^n - r^n)^{1/n}$ gives $s^{n-1}ds = -r^{n-1}dr$.}
\end{rem}

  \section{Symmetric kernels}

   The kernel defined by \eqref{modf_kernel1} is not symmetric in general,  even if $K(x,y)$ is symmetric. {It is symmetric if} the identity $K(\varphi_k(x), y) = K(x, \varphi_k(y))$ holds, for all $x,y \in \Omega $ and  $k=1,2, \cdots, {m}$.

    We now consider some special cases, where it is possible to obtain this property.

    Suppose    $J: {[0,\infty)} \to {[0,\infty)}$ {is such that $J(\|\cdot\|) \in L^1(\R^n)$, $\int_{\R^n} J(\|z\|)\,dz = 1$ and $J(s) = 0$ for $s \geq \delta$}, and $K(x,y) = J(\|x-y\|)$ is of convolution type. {Then $K$ satisfies \eqref{outflux} and $K(x,y) = 0$ if $\|x-y\|\geq\delta$.} Define the modified kernel as in \eqref{modf_kernel1}, {with $U_k = \Omega$}, which now reads
   \begin{equation} \label{modf_kernel2}  \tilde{K}(x,y)  = J(\|x-y\|) +  J(\| \varphi_1(x)- y \|) + \cdots  +  J(\| \varphi_{{m}}(x)- y \|).
    \end{equation}

    Then, one can prove the following result:
  \begin{lemma} \label{modf_kernel_isometry}
   Let $\Omega \subset \R^n$ {be a bounded domain with $|\partial\Omega| = 0$} and suppose there exist maps
    ${\varphi_1}, \dots, {\varphi_{{m}}}: {\R^n \to \R^n}$ such that
 \begin{enumerate}
  \item $\varphi_k(\Omega ) \subset \Omega^c$, for $k=1,2, \cdots, {m}$, {and the sets $\varphi_k(\Omega)$ are pairwise disjoint;}
   \item {$\{ z \in \Omega^c : \dist(z, \Omega) < \delta \} \subset \bigcup_{k=1}^m \overline{\varphi_k(\Omega)}$;}
 \item  $ \varphi_k$, $k=1,2, \cdots, {m}$,   are {involutive} isometries, i.e., $\varphi_k \circ \varphi_k = \textrm{id}$ and $\|\varphi_k(x) - \varphi_k(y)\| = \|x-y\| $ {for all $x, y \in \R^n$}.
  \end{enumerate}
  Then,  the kernel defined by {\eqref{modf_kernel2}}
  is  symmetric and
  \[ \int_{\Omega}   \tilde{K}(x,y)\, d \,x  = {\int_{\Omega}   \tilde{K}(y,x)\, d \,x =}    \int_{\R^n}   K(x,y)\, d \,x =1,\]
  for any $y \in  \Omega$. {In particular, $\tilde K$ satisfies \eqref{outflux_bounded} and \eqref{hyp_2}, and the constants are equilibria of the corresponding Neumann problem.}
     \end{lemma}
  \proof {Isometries of $\R^n$ are affine maps with $|\det D\varphi_k| = 1$; hence they preserve Lebesgue measure and $\overline{\varphi_k(\Omega)} \setminus \varphi_k(\Omega) \subset \varphi_k(\partial\Omega)$ has measure zero. Thus the hypotheses of Lemma \ref{const_integral} hold with $U_k = \Omega$, which gives the first identity once symmetry is proved.} We have
  \begin{align*}
  \tilde{K}(x,y) & = J(\|x-y\|) +  J(\| \varphi_1(x)- y \|) + \cdots  +  J(\| \varphi_{{m}}(x)- y \|)  \\
   & = J(\|x-y\|) +  J(\| \varphi_1^2(x)- \varphi_1(y) \|) + \cdots  +  J(\| \varphi_{{m}}^2(x)- \varphi_{{m}}(y) \|)  \\
   & = J(\|x-y\|) +  J(\| x - \varphi_1(y) \|) + \cdots  +  J(\| x- \varphi_{{m}}(y) \|)  \\
   & = \tilde{K}(y,x).
  \end{align*}
  \eproof

  The conditions of {Lemma \ref{modf_kernel_isometry}} can be met at least in the case of {rectangles and rectangular boxes, and of the regular polygons that tile the plane (the equilateral triangle, the square and the regular hexagon); in $\R^3$, the cube is the only regular polyhedron that tiles space.} The maps $\varphi_k$ can be defined by reflections about the lines (or planes) supporting the sides (or faces) and {by point reflections or reflections about} lines (planes) through the vertices with suitable inclinations. The case of the {square} is {especially} simple, as can be seen in the following example.

     \begin{ex}
   In the square ${\Omega = \,]-1,1[ \times ]-1,1[}$, we may define
  \[
   {\tilde K(x,y) := J(\|x-y\|) + \sum_{k=1}^4 J(\|\varphi_k(x)-y\|) + \sum_{k=1}^4 J(\|\psi_k(x)-y\|)},
   \]
   where
   \begin{itemize}
   \item $ \varphi_1(y_1, y_2)=  (2-y_1, y_2) $,
  $ \varphi_2(y_1, y_2)= (-2 -y_1, y_2)  $,  $ \varphi_3(y_1, y_2)=
   (y_1, 2-y_2)$  and $ \varphi_4(y_1, y_2)= (y_1, -2-y_2)$ are the reflections with respect to the lines  $y_1 =1$, $y_1=-1$,  $y_2 = 1$ and $y_2 = -1$,
    respectively;
    \item  $ \psi_1(y_1, y_2) = (2-y_1, 2-y_2) $,
    $ \psi_2(y_1, y_2) = (-2-y_1, 2-y_2) $, $ \psi_3(y_1, y_2) = (-2-y_1, -2-y_2) $, $ \psi_4(y_1, y_2) = (2-y_1, -2-y_2) $
      are the reflections with respect to the corners $(1,1)$,
      $(-1,1)$, $(-1,-1) $ and  $(1,-1)$ {(note that $\psi_1 = \varphi_1\circ\varphi_3$, etc.). Alternatively, one may use the reflections with respect to the lines $y_1+y_2=2$, $y_2-y_1=2$, $y_1+y_2=-2$ and $y_1-y_2=2$, which are different maps, but also satisfy the hypotheses of Lemma \ref{modf_kernel_isometry}.}
      \end{itemize}
 It is easy to check that these maps satisfy the hypotheses of
  Lemma \ref{modf_kernel_isometry}, {provided $\delta \leq 2$, the side of the square. For a rectangle, the condition is that $\delta$ does not exceed the smaller side.}
    \end{ex}

\begin{figure}[!h]
 \centering
 \includegraphics[width= 0.9 \columnwidth]{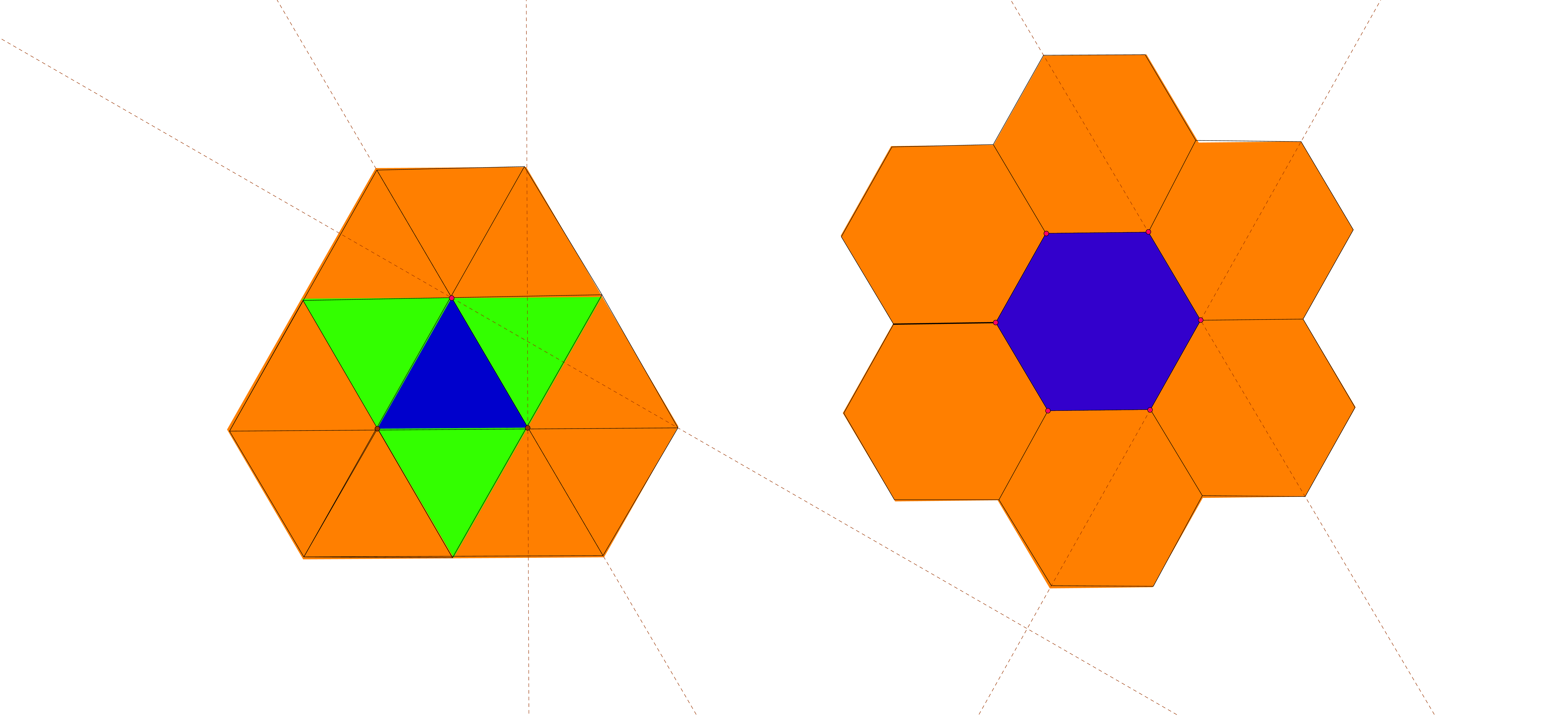}
 \caption{ Reflections for the triangle  and the hexagon, with some lines of reflection shown. }
 \label{triangle_hexagon}
\end{figure}

   {The cases of the regular hexagon and of the equilateral triangle, with side $\ell$, are similar, as illustrated in Figure \ref{triangle_hexagon}. For the hexagon, the six reflections about the lines supporting the sides suffice, and one may take $\delta \leq \ell$. For the triangle, besides the three reflections about the sides, at each vertex $v$ one needs three further maps, sending the triangle onto the remaining triangles of the tiling around $v$: the point reflection through $v$ and two reflections about lines through $v$ (the rotations by $\pm 2\pi/3$ about $v$ are not involutions, but, since the triangle is symmetric, the same image triangles are obtained from reflections). Here one must take $\delta \leq \frac{\sqrt3}{2}\ell$, the height of the triangle: indeed, if $M$ is the midpoint of a side and $D$ the vertex of the reflected triangle opposite to that side, then $D$ lies on the boundary of the union of the images and $\|M - D\| = \frac{\sqrt3}{2}\ell$.}

 {One obvious  question we have not been able to answer here is whether a kernel {satisfying} the {conclusion} of Lemma \ref{modf_kernel_isometry} exists in  more general domains.  Also, other boundary conditions, like   nonhomogeneous Neumann or Robin conditions can be considered (see \cite{ana}, \cite{costica}, for instance)}. {Another natural question is whether, after the usual rescaling $K_\epsilon(x,y) = \epsilon^{-n-2}K(x/\epsilon,y/\epsilon)$, the solutions of the Neumann problem with the reflected kernel converge to those of the local heat equation with Neumann boundary condition, as is known for the formulation with both integrals restricted to $\Omega$ \cite{cerw}.}

  \section*{Acknowledgement}
  The author {wishes} to thank the anonymous referee for her/his remarks and suggestions. 
  
\noindent   During the revision of this work, the author used the AI assistant Claude Opus 5.5 (Anthropic) to check the arguments and suggest corrections and simplifications, in particular in the proofs of Section 2. All suggestions were verified by the author, who takes full responsibility for the final text.


\begin{thebibliography}{99}

\bibitem{ana} {Ana-Maria Moșneagu, \emph{On some local and nonlocal reaction-diffusion models with Robin boundary conditions}, Discrete and Continuous Dynamical Systems - S, 2023, 16(1): 104-125. doi: 10.3934/dcdss.2022161.}

\bibitem{belletini} G. {Bellettini}, A. De Masi, E. Presutti, \emph{Energy levels of a nonlocal functional}, {J. Math. Phys. 46 (2005), no.~8, 083302}, DOI:10.1063/1.1990107.

\bibitem{bernal} A. {Rodr\'{\i}guez-Bernal}, S. {Sastre-G\'omez},  \emph{Nonlinear {nonlocal reaction-diffusion problem with local reaction}}, Discrete {Contin.} Dyn. Syst. {42} {(2022)}, no.~4, 1731--1765,
doi:10.3934/dcds.2021170.

\bibitem{cerw} {C. Cort\'azar, M. Elgueta, J. D. Rossi, N. Wolanski, \emph{How to approximate the heat equation with Neumann boundary conditions by nonlocal diffusion problems}, Arch. Ration. Mech. Anal. 187 (2008), 137--156.}

\bibitem{costica} {Costică Moroşanu, Bianca Satco, \emph{Qualitative and quantitative analysis for a nonlocal and nonlinear reaction-diffusion problem with in-homogeneous Neumann boundary conditions}, Discrete and Continuous Dynamical Systems - S, 2023, 16(1): 1-15. doi: 10.3934/dcdss.2022042.}

 \bibitem{fife} P. Fife, \emph{Some nonclassical trends in parabolic and parabolic-like evolutions}, {in: Trends in Nonlinear Analysis}, pp. 153--191, Springer, Berlin, 2003.

\bibitem{henry} D. B.  Henry, \emph{Perturbation of the {Boundary in Boundary-Value Problems} of Partial Differential Equations}, London Mathematical Society Lecture Note Series {318}, Cambridge University Press, {Cambridge, 2005}.

\bibitem{hutson} V. Hutson, S. Martinez, K. Mischaikow, G. T. Vickers, \emph{The evolution of dispersal}, J. Math. Biol. 47 {(2003)}, 483--517, DOI: 10.1007/s00285-003-0210-1.

\bibitem{rossi}  {F. Andreu-Vaillo, J. M. Maz\'on, J. D. Rossi, J. J. Toledo-Melero}, \emph{Nonlocal Diffusion Problems}, {Mathematical Surveys and Monographs 165}, American Math. Society, Real Sociedad Matem\'atica {Espa\~nola}, 2010.

\end{thebibliography}
\end{document}